\documentclass[12pt]{article}
\usepackage{graphics}
\newtheorem{THEOREM}{Theorem}[section]
\newtheorem{COROLLARY}{Corollary}[section]
\newcommand{\qed}{\ \rule[-1pt]{4pt}{8pt} 

                                         \vspace{2ex} }
\newenvironment{PROOF}{

                       \noindent{\bf Proof}.}{\qed}
\newcounter{labelflag} 

\newcommand{\Label}[1]{
                       \ifnum\thelabelflag=1 
                          \ifmmode  
                             \makebox[0in][l]{\qquad\fbox{\rm#1}}
                          \else
                             \marginpar{\vspace{0.7\baselineskip}
                                        \hspace{-1.1\textwidth}
                                        \fbox{\rm#1}}
                          \fi 
                       \fi
                       \label{#1} 
                      }
\newcommand{\be}{\begin{equation}}
\newcommand{\ee}{\end{equation}}
\newcommand\evec{\mbox{\boldmath{$e$}}}
\newcommand\nvec{\mbox{\boldmath{$n$}}}
\newcommand\xvec{\mbox{\boldmath{$x$}}}
\newcommand\Avec{\mbox{\boldmath{$A$}}}
\newcommand\Bvec{\mbox{\boldmath{$B$}}}
\newcommand\Evec{\mbox{\boldmath{$E$}}}
\newcommand\Fvec{\mbox{\boldmath{$F$}}}
\newcommand\Hvec{\mbox{\boldmath{$H$}}}
\newcommand\Jvec{\mbox{\boldmath{$J$}}}
\newcommand\Mvec{\mbox{\boldmath{$M$}}}
\newcommand\Uvec{\mbox{\boldmath{$U$}}}
\newcommand{\eps}{\varepsilon}
\newcommand{\grad}{\nabla}
\newcommand{\curl}{\nabla\times}
\renewcommand{\div}{\nabla\cdot}
\newcommand{\half}{\mathrm{\frac{1}{2}}}
\begin{document}


\begin{center}
\begin{large}
\textbf{The frozen-field approximation and the} \\[1ex]
\textbf{Ginzburg--Landau equations of superconductivity}
\end{large}

Hans G.\ Kaper\footnote{
Mathematics and Computer Science Division,
Argonne National Laboratory,
Argonne, IL 60439, USA
(\texttt{kaper@mcs.anl.gov})
}
and Henrik Nordborg\footnote{
James Franck Institute,
The University of Chicago,
5640 South Ellis Avenue,
Chicago, IL 60637, USA
(\texttt{Henrik\_Nordborg@anl.gov})
}
\end{center}

\paragraph{Abstract.}
The Ginzburg--Landau (GL) equations of superconductivity
provide a computational model for the study of magnetic
flux vortices in type-II superconductors.
In this article we show through numerical examples
and rigorous mathematical analysis that the GL model
reduces to the frozen-field model
when the charge of the Cooper pairs
(the superconducting charge carriers)
goes to zero while the applied field
stays near the upper critical field.


\section{Introduction\label{s-intro}} 
\setcounter{equation}{0}
Superconducting materials hold great promise for
technological applications.
Especially since the discovery of the so-called
high-temperature superconductors in the 1980s,
much research has been devoted to understanding
the behavior of these new materials.
While conventional superconductors require
liquid helium (3--4 degrees Kelvin) to remain
in the superconducting state,
high-temperature superconductors can be cooled
with liquid nitrogen (76 degrees Kelvin)---a
clear economic advantage.
Unfortunately, high-temperature superconductors
are ceramic materials, which are difficult to
manufacture into films and wires, but progress
is being made all the time.

High-temperature superconductors belong to the class
of type-II superconductors.
Unlike type-I superconductors, type-II superconductors
can sustain a magnetic flux in their interior,
but this flux is restricted to quantized
amounts---filaments that are encircled by a current.
The current shields the magnetic flux from
the bulk, which is perfectly superconducting.
The configuration resembles that of a vortex
in a fluid, and the superconductor is said to be
in the \textit{vortex state}.

\begin{figure}
\begin{center}
\resizebox{3.0in}{!}{\mbox{\includegraphics{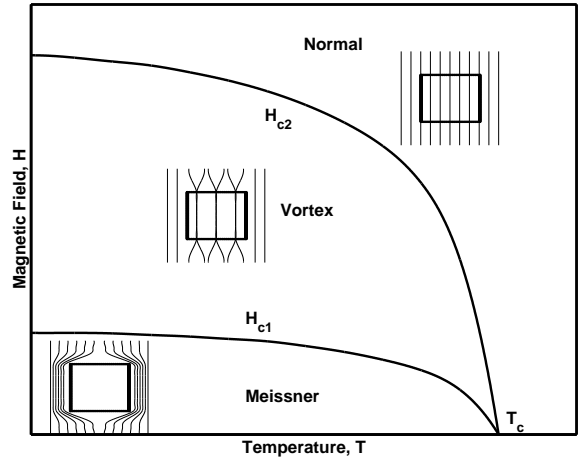}}}
\caption{Phase diagram of a type-II superconductor.
         \label{f-sctype2}}
\end{center}
\end{figure}
Figure~\ref{f-sctype2} gives a sketch of
the phase diagram of a type-II superconductor
in the neighborhood of $T_c$,
the \textit{critical temperature}.
The two-dimensional phase space is
spanned by the temperature $T$ and
the (magnitude of the) magnetic field $H$
and is roughly divided into three subregions.
Each subregion corresponds to a particular state:
the perfectly superconducting (Meissner) state
below the lower critical field $H_{c1}$,
where no magnetic field can penetrate the medium;
the normal state above the upper critical field $H_{c2}$,
where the superconductor behaves like a normal metal;
and the intermediate vortex state.
Above the critical temperature $T_c$
all superconducting properties are lost.

The vortices, and especially their dynamics,
determine the current-carrying capabilities
of a superconductor.
Much effort, both experimental and theoretical,
is therefore being spent on the study of vortex dynamics
and, especially, mechanisms to inhibit vortex
motion when the superconductor is subject to
currents and fields.
By ``pinning'' the vortices, one prevents energy
dissipation and, hence, loss of superconductivity.

\begin{figure}
\begin{center}
\resizebox{2.0in}{!}{\mbox{\includegraphics{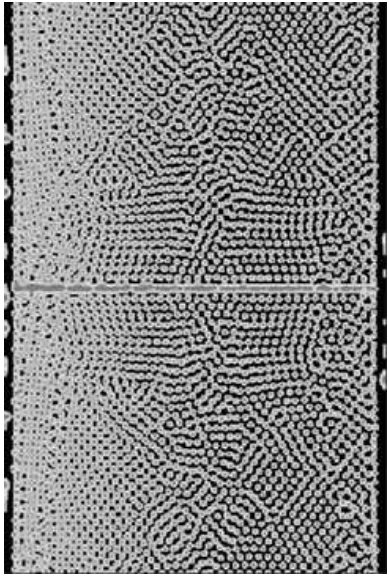}}}
\caption{Vortex configuration in two dimensions.
         \label{f-2d}}
\end{center}
\end{figure}

\begin{figure}
\begin{center}
\resizebox{1.0in}{!}{\mbox{\includegraphics{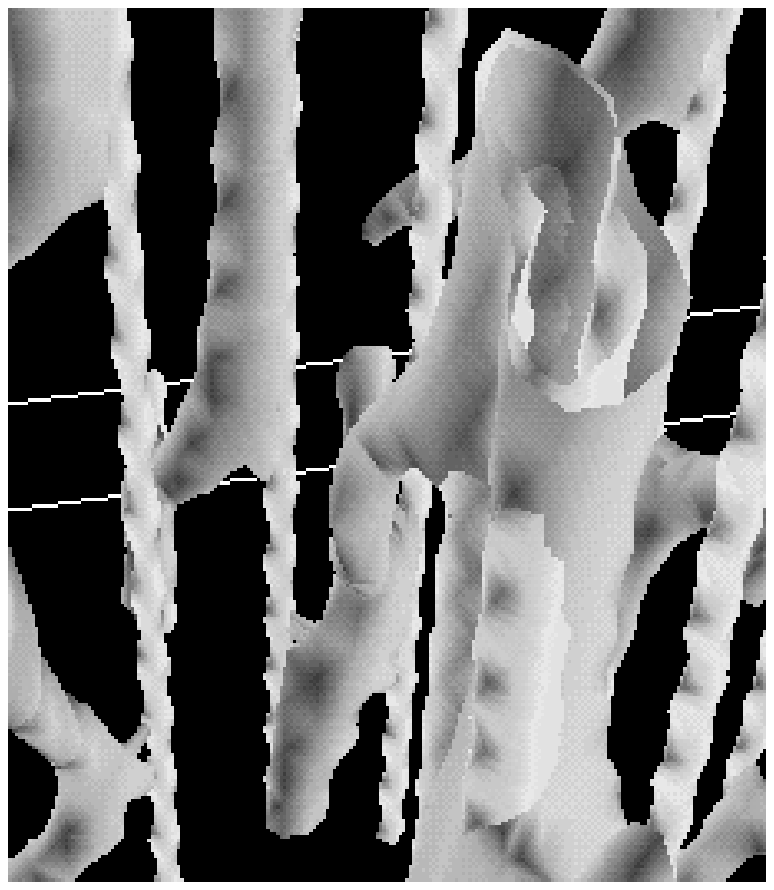}}}
\hspace*{1em}
\resizebox{1.0in}{!}{\mbox{\includegraphics{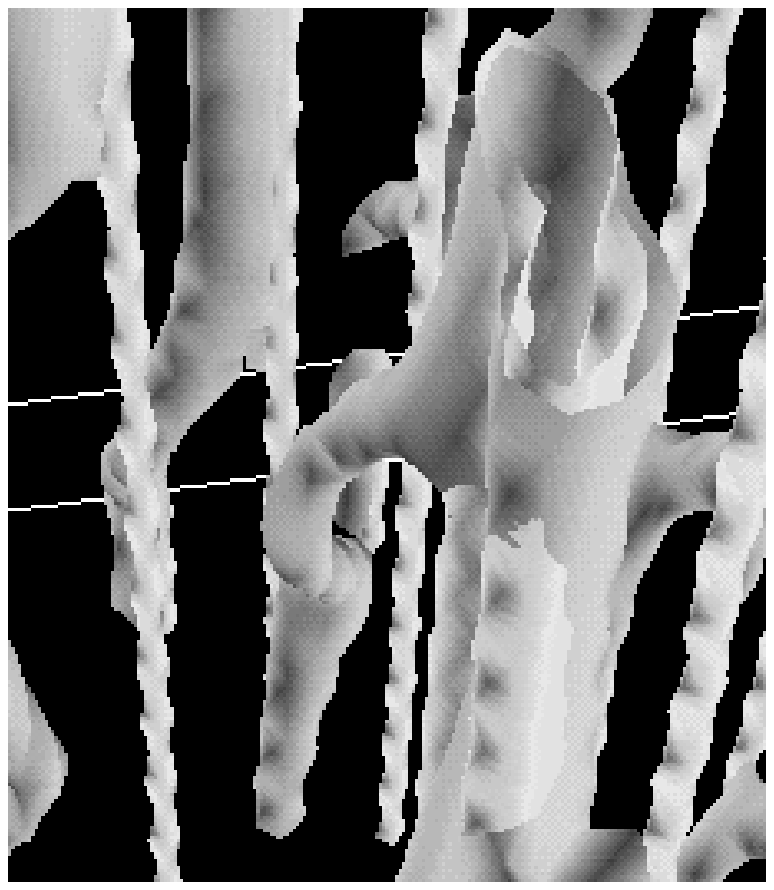}}}
\hspace*{1em}
\resizebox{1.0in}{!}{\mbox{\includegraphics{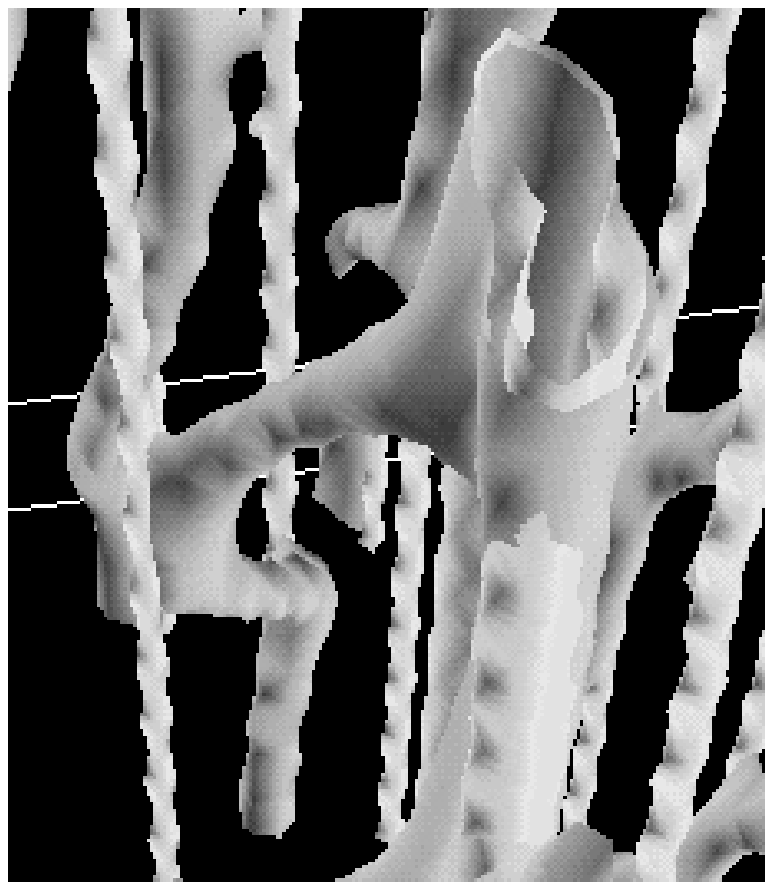}}} \\
\vspace{2ex}
\resizebox{1.0in}{!}{\mbox{\includegraphics{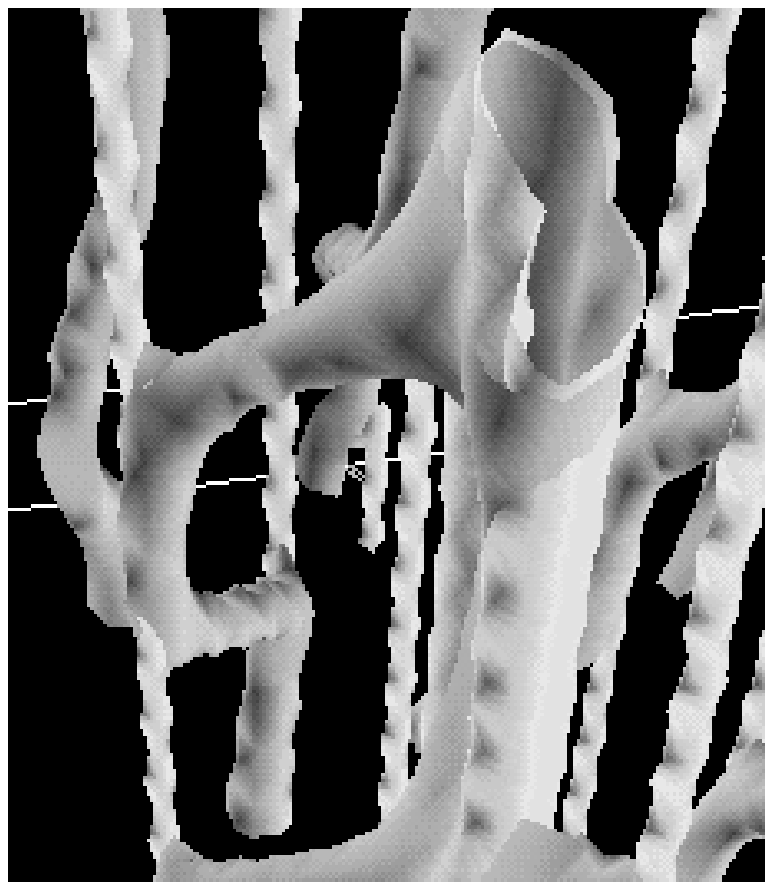}}}
\hspace*{1em}
\resizebox{1.0in}{!}{\mbox{\includegraphics{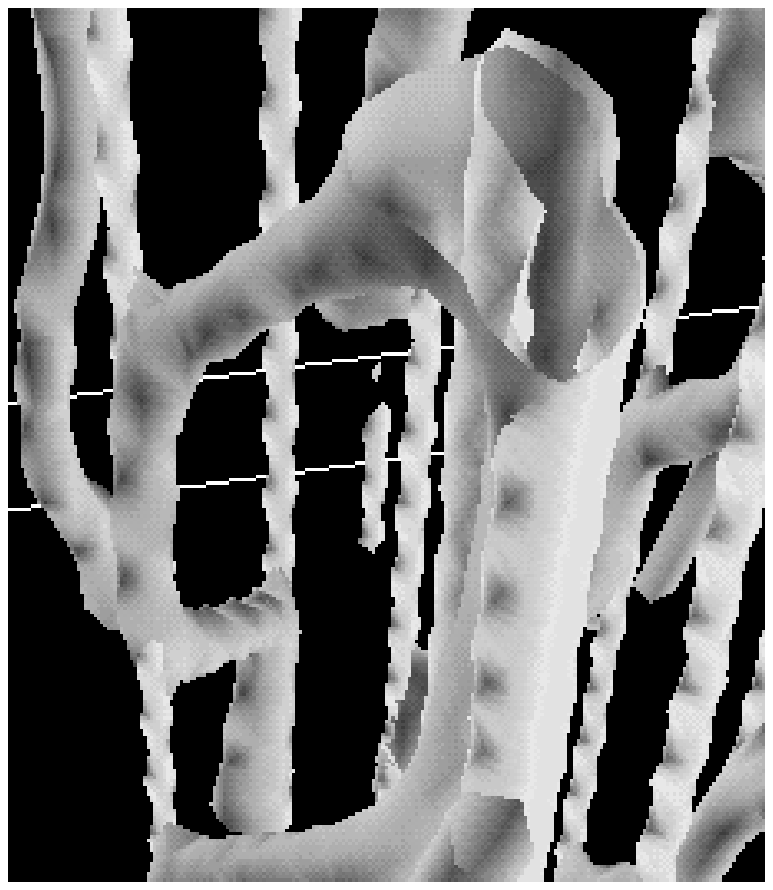}}}
\hspace*{1em}
\resizebox{1.0in}{!}{\mbox{\includegraphics{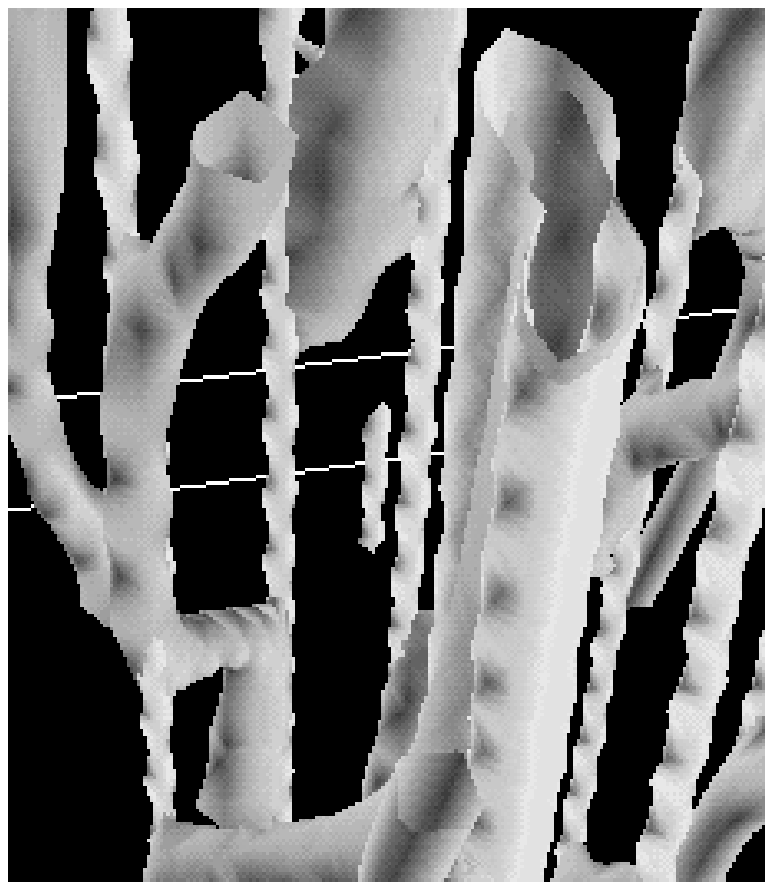}}}
\caption{Kinking-induced motion of vortices through splayed columnar defects.
         \label{f-3d}}
\end{center}
\end{figure}

Vortices can be studied computationally
at various levels of detail using different models.
The Ginzburg--Landau (GL) model gives
a field (continuum) description that,
although phenomenological and not based on
any microscopic quantum-mechanical theory,
has been used successfully to study both
the dynamics and the structure of vortex systems
in realistic superconductor configurations~\cite{PRL,PRB}.
Figures~\ref{f-2d} and \ref{f-3d} give two examples
of computational results obtained with the GL equations.
They illustrate both the effectiveness and
the difficulties of such calculations.

Figure~\ref{f-2d} shows a vortex configuration in
a two-dimensional cross section of a twinned
superconducting crystal, which was computed
from a steady-state solution of the GL equations.
The twin boundary (an irregularity in the structure
of the crystal) is visible in the horizontal line
through the center; it acts as a pinning site for the vortices.
The field is perpendicular to the plane of the cross section,
which measures $128 \times 192$ coherence lengths
(a characteristic length of the order of microns).
Each dot corresponds to a vortex intersecting the plane
of the cross section; the entire configuration has
approximately 2,700 vortices.
The figure shows the level of detail one can achieve
with the GL model, given sufficient computing power.
At the same time, it illustrates the level of
computational complexity one faces if one uses
the GL model.

Figure~\ref{f-3d} shows a series of snapshots of
a vortex configuration in three dimensions,
also computed with the GL model.
The objective of this computation was to simulate
vortex motion through columnar defects and study
the potential of the latter as pinning sites.
The defects are visible as twisted straight lines.
The vortices are the flexible tube-like structures;
they move from one defect to another under the influence
of external forces.
The figure shows the motion of a vortex that is
originally pinned on a defect.
The vortex develops a loop, the loop peels off,
the loop expands in both directions in a
traveling-wave-like scenario, and gradually
the entire vortex transfers to the
next available defect.

Numerical simulations provide the only way to study
vortex dynamics at this level of detail.
They are an invaluable tool for fundamental research,
complementing experiment and theory.
Numerical simulations of realistic superconductors
based on the GL model, like the ones illustrated
in Figs.~\ref{f-2d} and \ref{f-3d} are, however,
extremely time consuming, and it is desirable
to use simpler models whenever possible.
Here, we focus on the ``frozen-field model,''
which is still a continuum model and the closest
approximation to the full GL model.
In the frozen-field model, the superconducting
phenomena are decoupled from the electromagnetic field,
and the latter is prescribed through a vector potential.
The frozen-field model is much simpler and has been used
successfully for numerical simulations of vortex
systems~\cite{Igor}.

In this article, we prove that the frozen-field model
is obtained as the asymptotic limit
of the GL model when the charge of the Cooper pairs
(the superconducting charge carriers)
goes to zero while the applied magnetic field
stays near the upper critical field.
(The upper critical field itself depends
on the charge of the Cooper pairs and
increases as the latter decreases.)
Because the temperature is constant in the GL model,
this limit corresponds to fixing the temperature $T$ 
and moving up vertically through the vortex regime
to the curve labeled $H_{c2}$ in the phase diagram
of Fig.~\ref{f-sctype2}.
The convergence rate is second order in the small parameter.

For more background on the physics of superconductivity
we refer the reader to the monograph by Tinkham~\cite{Tinkham}.
The original source for the GL equations
of superconductivity is~\cite{Ginzburg-Landau}.
A good introduction to the mathematics of the
GL equations is~\cite{Du-92}.
The dynamics of the GL equations
have been studied by several authors;
see~\cite{Du-94,Tang-Wang,Fleck-Kaper-Takac-98}
and the references cited therein.
The present investigation is closely related to
the work of Du and Gray~\cite{Du-Gray-96}.

Section~\ref{s-gl} introduces
the Ginzburg--Landau equations,
Section~\ref{s-num} contains the numerical results
and Section~\ref{s-aa} the analysis.

\section{The Ginzburg--Landau equations\label{s-gl}} 
\setcounter{equation}{0}
In the Ginzburg--Landau theory of superconductivity,
the state of a superconducting medium is described by
a complex scalar-valued \textit{order parameter}~$\psi$ and
a real vector-valued \textit{vector potential}~$\Avec$.
If the state varies with time, a third variable---the
\textit{electric potential}~$\phi$---is necessary
to fully describe the electromagnetic field.
The evolution of the state variables is governed
by the time-dependent Ginzburg--Landau (TDGL) equations,
\be
  \gamma \hbar
  \left(
  \frac{\partial}{\partial t}
  +
  \mathrm{i} \frac{q_s}{\hbar} \phi
  \right) \psi
  +
  \frac{1}{2m_s}
  \left(
  \frac{\hbar}{\mathrm{i}} \grad - \frac{q_s}{c} \Avec
  \right)^2 \psi
  + \alpha \psi + \beta | \psi |^2 \psi
  = 0 ,
  \Label{TDGL-p}
\ee
\be
  \sigma
  \left(
  -\frac{1}{c} \frac{\partial \Avec}{\partial t}
  - \grad \phi
  \right)
  - \frac{c}{4\pi} \curl\curl\Avec
  +
  \Jvec_s
  + \frac{c}{4\pi} \curl\Hvec
  = 0 ,
  \Label{TDGL-A}
\ee
where the supercurrent density $\Jvec_s$ is a nonlinear function
of $\psi$ and $\Avec$,
\be
  \Jvec_s =
  \frac{q_s \hbar}{2 \mathrm{i} m_s}
  (\psi^* \grad \psi - \psi \grad \psi^*)
  -
  \frac{q_s^2}{m_s c} | \psi|^2 \Avec
  = \frac{q_s}{m_s}
  \Re \left[ \psi^* \left(
  \frac{\hbar}{\mathrm{i}} \grad - \frac{q_s}{c} \Avec
  \right) \psi \right] .
  \Label{GL-Js}
\ee
These equations are supplemented by the boundary conditions,
\be
  \nvec \cdot \Jvec_s = 0 , \quad
  \nvec \times (\curl \Avec) = \nvec \times \Hvec .
  \Label{GL-bc}
\ee
Here,
$\Hvec$ is the \textit{applied magnetic field},
which we assume to be time independent.
The constants $m_s$ and $q_s$ are
the mass and charge, respectively,
of a Cooper pair
(the superconducting charge carriers,
also referred to as superelectrons);
$c$ is the speed of light;
and $\hbar$ is Planck's constant divided by $2\pi$.
A Cooper pair is made up of two electrons,
each with charge $-e$ ($e$~is the elementary charge);
hence, $q_s$ is negative, $q_s = -2e$.

The parameters $\alpha$ and $\beta$ are material parameters;
$\alpha$ changes sign at the critical temperature~$T_c$,
$\alpha (T) < 0$ for $T < T_c$ (superconducting state)
and
$\alpha (T) > 0$ for $T > T_c$ (normal state);
$\beta$ is only weakly temperature dependent
and positive for all $T$.
The remaining parameters are
$\sigma$, the normal state conductivity,
and $\gamma$, the mobility coefficient.
The latter is dimensionless and related to
the diffusion coefficient~$D$,
$\gamma = \hbar/2m_sD$.

The boundary conditions~(\ref{GL-bc}) express the fact
that superelectrons cannot leave the superconductor.
Also, if no surface currents are present,
the tangential components of the magnetic field
must be continuous across the boundary.

The parameters $\alpha$ and $\beta$
are defined phenomenologically,
but they can be expressed in terms of measurable quantities,
such as the superconducting \textit{coherence length}~$\xi$
and the London \textit{penetration depth}~$\lambda$,
\be
  \xi = \left( \frac{\hbar^2}{2m_s |\alpha|} \right)^{1/2} , \quad
  \lambda = \left( \frac{m_s c^2 \beta}{4\pi q_s^2 |\alpha|} \right)^{1/2} .
  \Label{xi-lambda}
\ee
The coherence length and the London penetration depth
define the respective characteristic length scales for
the order parameter and the magnetic induction.
Both depend on the temperature~$T$ and diverge as $T$
approaches the critical temperature $T_c$, because
of the factor $|\alpha|^{-1/2}$.
However, their ratio is, to a good approximation,
independent of temperature.
This ratio is the \textit{Ginzburg-Landau parameter},
\be
  \kappa = \lambda / \xi .
  \Label{GL-parm}
\ee
In high-$T_c$ superconductors, $\kappa$ is of
the order of 50--100.

The electromagnetic variables are the
\textit{magnetic induction}~$\Bvec$,
the \textit{current density}~$\Jvec$, and
the \textit{electric field}~$\Evec$;
they are given in terms of $\Avec$ and $\phi$
by the expressions
\be
  \Bvec = \curl \Avec , \quad
  \Jvec = \frac{c}{4\pi} \curl \curl \Avec , \quad
  \Evec = - \frac{1}{c} \frac{\partial \Avec}{\partial t} - \grad \phi .
  \Label{TDGL-BJE}
\ee
Equation~(\ref{TDGL-A}) is essentially Amp\`{e}re's law,
$\Jvec = (c/4\pi) \curl \Bvec$,
where the current $\Jvec$
is the sum of the supercurrent $\Jvec_s$,
the transport current~$\Jvec_t = (c/4\pi) \curl \Hvec$,
and a ``normal'' current~$\Jvec_n = \sigma \Evec$
(Ohm's law).
Hence, the GL model uses a quasistatic version
of Maxwell's equations, where the time derivative
of the electric field is ignored.

The TDGL equations were first given by
Schmid~\cite{Schmid-66} in 1966
and subsequently derived from
the microscopic theory of superconductivity
by Gor'kov and Eliashberg~\cite{Gorkov-Eliashberg-68}.
Our notation is the same as in
Gor'kov and Kopnin~\cite{Gorkov-Kopnin-76}.

The solution of the TDGL equations is not unique.
Any solution $(\psi, \Avec, \phi)$ defines
a family of solutions $G_\chi (\psi, \Avec, \phi)$
parameterized by a sufficiently smooth function
$\chi$ of space and time,
\be
  G_\chi: (\psi, \Avec, \phi) \mapsto
     \left( \psi \textrm{e}^{i(q_s/\hbar c) \chi},
            \Avec + \grad \chi,
            \phi - \frac{1}{c} \frac{\partial \chi}{\partial t} \right) .
  \Label{G-gauge}
\ee
This property is known as \textit{gauge invariance};
the function $\chi$ is known as a \textit{gauge function}.
Gauge invariance does not affect the
physically measurable quantities
(the magnetic induction~$\Bvec$,
the magnetization~$\Mvec = \Bvec - \Hvec$,
and the current density~$\Jvec$).
Uniqueness requires an additional constraint,
which is imposed through a gauge choice.
The choice of a proper gauge for the TDGL equations
has been a subject of considerable debate.
The choice is a matter of convenience
and may depend on the problem under investigation.
In this article we adopt a gauge in which,
at any time, the electric potential and
the divergence of the vector potential
satisfy the identity
\be
  \sigma \phi + (c/4\pi) \div \Avec = 0
  \Label{TDGL-gauge}
\ee
everywhere in the domain, while $\Avec$ is
tangential at the boundary.
This choice is realized by identifying the gauge $\chi$
with a solution of the linear parabolic equation
\be
  \frac{\sigma}{c} \frac{\partial \chi}{\partial t}
  - \frac{c}{4\pi} \Delta \chi
  = \sigma \phi + \frac{c}{4\pi} \div \Avec ,
  \Label{chi-eqn}
\ee
subject to the condition
$\nvec \cdot \grad \chi = - \nvec \cdot \Avec$
on the boundary.
In~\cite{Fleck-Kaper-Takac-98}, it was shown that the
TDGL equations, subject to the constraint~(\ref{TDGL-gauge}),
define a dynamical system under suitable regularity
conditions on $\Hvec$.
(In the more general case, where $\Hvec$ varies not
only in space but also in time, the TDGL equations
define a dynamical process.)
This dynamical system has a global attractor,
which consists of the stationary points of the
dynamical system and the heteroclinic orbits
connecting such stationary points.
Furthermore, it was shown that every solution
on the attractor satisfies the condition $\div\Avec = 0$
(and, therefore, also $\phi = 0$).
Thus, in the limit as $t \to \infty$,
every solution of the TDGL equations satisfies
the GL equations in the London gauge.

\subsection{Nondimensional TDGL equations\label{ss-nondim-tdgl}}
In this section, we render the TDGL equations dimensionless
by choosing units for the independent and dependent
variables.
Since we are interested in the collective behavior
of vortices in the bulk of a superconductor in
the limit of weak coupling ($q_s \to 0$),
we take care to choose the units in such a way
that they remain of order one as $q_s \to 0$.
(We recall that $q_s$ is negative, $q_s = - 2e$.)

As $q_s \to 0$, the coherence length~$\xi$
remains of order one,
while the penetration depth~$\lambda$ increases
like $|q_s|^{-1}$;
see Eq.~(\ref{xi-lambda}).
This suggests taking the coherence length
$\xi$ as the unit of length.

To maintain the diffusion coefficient
$D = \hbar/2\gamma m_s = \xi^2 (\gamma\hbar/|\alpha|)^{-1}$
at order one,
we measure time in units of $\gamma \hbar/|\alpha|$.

The real and imaginary parts of the order parameter
are conveniently measured in units of
$\psi_0 = (|\alpha|/\beta)^{1/2}$,
which is the value of $\psi$ that minimizes
the free energy in the absence of a field.

Next, consider the magnetic field.
A fundamental quantity in the theory of type-II
superconductors is the flux quantum~$\Phi_0$,
\be
  \Phi_0 = \frac{hc}{|q_s|} = 2\pi \frac{\hbar c}{|q_s|} .
\ee
The flux quantum is the unit of magnetic flux carried by a vortex.
Together with the coherence length and penetration depth,
it defines three characteristic field strengths:
the \textit{lower critical field}~$H_{c1}$,
the \textit{thermodynamical critical field}~$H_c$, and
the \textit{upper critical field}~$H_{c2}$,
\be
  H_{c1} = \frac{\Phi_0}{4\pi \lambda^2 \ln \kappa} , \quad
  H_c = \frac{\Phi_0}{2\pi \xi \lambda \surd{2}} , \quad
  H_{c2} = \frac{\Phi_0}{2\pi \xi^2} .
\ee
Below~$H_{c1}$, a superconductor is in the
ideal superconducting (Meissner) state,
where it does not support magnetic flux in the bulk;
above~$H_{c2}$, it is in the normal state,
where the magnetic flux is distributed uniformly
in the bulk;
between~$H_{c1}$ and $H_{c2}$, it is in the vortex state,
where magnetic flux is quantized in vortex-like configurations
(see Fig.~\ref{f-sctype2}).
The thermodynamical critical field~$H_c$
is intermediate between $H_{c1}$ and $H_{c2}$
and is defined by the identity
$H_c^2/8\pi = \half \psi_0^2 |\alpha|$;
$H^2/8\pi$ is the energy per unit volume
associated with a field $H$, and
$\half \psi_0^2 |\alpha|$ is the minimum
condensation energy, which is attained
when $\psi = \psi_0$,
so these two quantities are in balance
when $H = H_c$.

As $q_s \to 0$, $H_{c1}$ goes to 0 like $|q_s|$,
$H_c$ remains of order one, and $H_{c2}$ grows
like $|q_s|^{-1}$.
This suggests that we define field strengths
in terms of $H_c$.
In fact, it is convenient to absorb
a factor $\surd 2$, so we adopt
$H_{c} \surd 2$ or, equivalently,
$\hbar c / \xi \lambda |q_s|$
as the unit for the magnetic field strength.

With the coherence length as the unit of length
and $H_c\surd 2$ as the unit of field strength,
it follows that the vector potential is measured
in units of $\xi H_c \surd 2$.
Furthermore, energy densities are measured
in units of $H_{c}^2/4\pi$, which is the same as
$|\alpha| \psi_0^2$.

Finally, we define the scalar potential $\phi$
in units of $(1/\gamma \psi_0^2\kappa |q_s|)(H_c^2/4\pi)$.
Notice that this unit remains of order one as $q_s \to 0$,
because $\kappa |q_s|$ is of order one.
On the other hand, the product $q_s \phi$,
which represents an energy density,
tends to zero as $q_s \to 0$.
(It remains finite on the scale of the penetration depth.)

Table~\ref{nondim} summarizes the relations
between the original variables and their
nondimensional (primed) counterparts.
We adopt the latter as the new variables
and work until further notice
on the nondimensional problem.
We omit all primes.
\begin{table}[ht]
\begin{center}
\begin{footnotesize}
\caption{Nondimensionalization.
  \label{nondim}}\vspace{1ex}
\begin{tabular}{||l|l||}\hline
Independent variables
  & $\xvec = \xi \xvec'$ \\
  & $t = (\gamma \hbar / |\alpha|) t'$ \\\hline
  & $\psi = \psi_0 \psi'$ \\
Dependent variables
  & $\Avec = (\xi H_c \surd 2) \Avec'$ \\
  & $\phi = (1/\gamma\psi_0^2\kappa|q_s|)(H_c^2/4\pi) \phi'$ \\\hline
  & $\Bvec = (H_c \surd 2) \Bvec'$ \\
Electromagnetic variables
  & $\Jvec = (c H_c \surd 2/ 4\pi \xi) \Jvec'$ \\
  & $\Evec = (1/\gamma \psi_0^2\kappa |q_s|)(H_c^2/4\pi \xi) \Evec'$ \\\hline
Applied field
  & $\Hvec = (H_c \surd 2) \Hvec'$ \\
Normal conductivity
  & $\sigma = (\gamma m_s c^2/2\pi \hbar) \sigma'$ \\\hline
\end{tabular}
\end{footnotesize}
\end{center}
\end{table}

The nondimensional TDGL equations are
\be
  \left( \frac{\partial}{\partial t} - \frac{\mathrm{i}}{\kappa} \phi \right) \psi
  - \left( \grad + \frac{\mathrm{i}}{\kappa} \Avec \right)^2 \psi
  - (1 - | \psi |^2) \psi
  = 0 ,
  \Label{tdgl-p}
\ee
\be
  \sigma \frac{\partial \Avec}{\partial t}
  - \Delta \Avec
  - \frac{1}{\kappa} \Jvec_s
  - \curl \Hvec
  = {\bf 0} ,
  \Label{tdgl-A}
\ee
where
\be
  \Jvec_s
  =
  - \frac{1}{2\mathrm{i}}
  (\psi^* \grad \psi - \psi \grad \psi^*)
  -
  \frac{1}{\kappa}
  |\psi|^2 \Avec
  =
  - \Im
  \left[ \psi^* \left( \grad + \frac{\mathrm{i}}{\kappa} \Avec \right) \psi \right] ,
  \Label{tdgl-Js}
\ee
with the corresponding gauge condition,
\be
  \sigma \phi +  \div \Avec = 0 .
 \Label{tdgl-gauge}
\ee
In deriving Eq.~(\ref{tdgl-A}), we have made use of
the gauge condition~(\ref{tdgl-gauge})
and the vector identity
\be
  \Delta \Avec = - \curl \curl \Avec + \grad (\div\Avec) .
\ee
If $\Omega$ is the domain occupied by
the superconducting material
(measured in units of $\xi$),
then Eqs.~(\ref{tdgl-p})--(\ref{tdgl-gauge})
must be satisfied everywhere $\Omega$.
At the boundary $\partial \Omega$ of $\Omega$,
we have the conditions
\be
  \nvec \cdot \Jvec_s = 0 , \quad
  \nvec \times (\curl \Avec) = \nvec \times \Hvec , \quad
  \nvec \cdot \Avec = 0 .
  \Label{tdgl-bc}
\ee
Here, $\nvec$ is the local unit normal vector.

The electromagnetic variables are given by the expressions
\be
  \Bvec = \curl \Avec , \quad
  \Jvec = \curl\curl \Avec , \quad
  \Evec = - \partial_t \Avec - \grad \phi .
  \Label{tdgl-BJE}
\ee
The values of the lower and upper critical fields are
\be
  H_{c1} = (2 \kappa \ln \kappa)^{-1} , \quad
  H_{c2} = \kappa .
  \Label{Hc12}
\ee

\subsection{Link variables\label{ss-link}}
The combination $\grad + (\mathrm{i}/\kappa) \Avec$ plays
a fundamental role;
we refer to it as the $\Avec$-gradient and write
\be
  \grad_{\bf A} = \grad + \frac{\mathrm{i}}{\kappa} \Avec .
  \Label{Agrad}
\ee
The $\Avec$-gradient defines the $\Avec$-Laplacian
(or ``twisted Laplacian''),
\be
  \Delta_{\bf A}
  =
  \grad_{\bf A} \cdot \grad_{\bf A}
  =
  \left( \grad + \frac{\mathrm{i}}{\kappa} \Avec \right)^2 .
  \Label{A-Laplacian}
\ee
The relation between the $\Avec$-Laplacian
and the ordinary Laplacian is most easily
illustrated by means of the \textit{link variables},
\begin{eqnarray}
  U_x (x,y,z)
  &=& \exp \left( \frac{\mathrm{i}}{\kappa}
  \int^x A_x (\xi,y,z) \,{\rm d}\xi \right) , \nonumber \\
  U_y (x,y,z)
  &=& \exp \left( \frac{\mathrm{i}}{\kappa}
  \int^y A_y (x,\eta,z) \,{\rm d}\eta \right) , \\
  U_z (x,y,z)
  &=& \exp \left( \frac{\mathrm{i}}{\kappa}
  \int^z A_z (x,y,\zeta) \,{\rm d}\zeta \right) . \nonumber
  \Label{tdgl-U}
\end{eqnarray}
(We omit the argument $t$.)
The integrals are evaluated
with respect to an arbitrary reference point.
Each $U_\mu$ ($\mu = x,y,z$) is complex valued
and unimodular, $U^*_\mu = U^{-1}_\mu$.
The vectors $\Avec$ and $\Uvec$ may be used interchangeably.
With a slight abuse of notation, we have
\be
  \Uvec = {\rm e}^{(\mathrm{i}/\kappa) \int \Avec} ,
  \quad
  \grad_{\bf A}
  =
  U^* \grad U ,
  \quad
  \Delta_{\bf A}
  =
  U^* \Delta U .
  \Label{Uvec}
\ee

\section{Numerical solution\label{s-num}}
\setcounter{equation}{0}
A parallel code for solving Eqs.~(\ref{tdgl-p})--(\ref{tdgl-bc})
has been developed as part of a project for
large-scale simulations of vortex dynamics
in superconductors.
Details on these simulations and
on the code will be published elsewhere;
here, we give only a brief overview
of the numerical methods and the results
of numerical simulations showing the behavior
of the solution as $\kappa$ increases.

The algorithm uses finite differences
on a staggered grid,
making all approximations accurate
to second order in the mesh widths,
and an implicit method for the time integration,
making the algorithm (essentially) unconditionally stable.
The code, written in C++, has been designed for
a multiprocessing environment;
it uses MPI for message passing.

We restrict the discussion to rectangular
two-dimensional configurations that are
periodic in one direction and open in the other.
The configurations are assumed to be infinite
in the third, orthogonal direction, which is also
the direction of the applied magnetic field,
$\Hvec = (0, 0, H_z)$.

\subsection{Discretization\label{ss-discrete}}
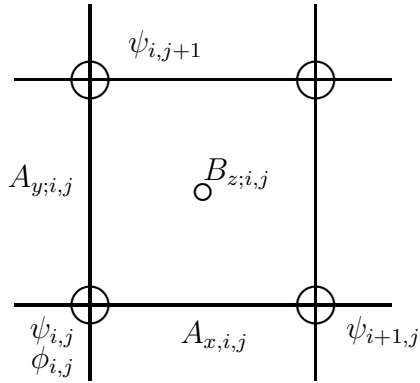
\begin{figure}[ht]
\begin{center}
\unitlength1cm
\begin{picture}(6,6)
\thicklines
\put(0,1){\line(1,0){5}}
\put(1,0){\line(0,1){5}}
\put(0,4){\line(1,0){5}}
\put(4,0){\line(0,1){5}}
\put(1,4){\circle{0.5}}
\put(4,1){\circle{0.5}}
\put(4,4){\circle{0.5}}
\put(1,1){\circle{0.5}}
\put(2.5,2.5){\circle{0.2}}
\put(0,0.4){\makebox(1,0.5){$\psi_{i,j}$}}
\put(0,0.0){\makebox(1,0.5){$\phi_{i,j}$}}
\put(2.7,2.5){\makebox(0.5,0.5){$B_{z;i,j}$}}
\put(2.2,0.5){\makebox{$A_{x,i,j}$}}
\put(0.1,2.4){\makebox(0.5,0.5){$A_{y;i,j}$}}
\put(4.4,0.4){\makebox(1,0.5){$\psi_{i+1,j}$}}
\put(1.5,4.2){\makebox(1,0.5){$\psi_{i,j+1}$}}
\end{picture}
\end{center}
\caption{Computational grid cell and definition of the discrete variables.
\label{f-cell}}
\end{figure}

\paragraph{Computational grid.}
The computational grid is uniform, with equal mesh sizes
in the $x$ and $y$ direction, $h_x = h_y = h$.
A vertex on the grid is denoted by $\xvec_{i,j} = (x_i, y_j)$
and is the point of reference for the grid cell
shown in Fig.~\ref{f-cell}.
The indices run through the values
$i=1, \ldots\,, N_x$ and $j=1, \ldots\,, N_y$.
We assume periodicity in the $x$ direction
and take the grid so the vertices with $j=1$ and $j=N_y$
are located on the open boundary of the superconductor.
Thus, the size of the domain is $S = N_x (N_y-1) h^2$.

\paragraph{Variables.}
The discrete variables are introduced so that
all derivatives are given by second-order accurate
central-difference approximations.
The scalar variables $\psi$ and $\phi$ are defined
on the vertices of the grid,
\be
  \psi_{i,j} = \psi(\xvec_{i,j}) , \quad
  \phi_{i,j} = \phi(\xvec_{i,j}) .
\ee
(We use the same symbol for the original field
and its discrete counterpart.)
Vectors are defined at the midpoints of the links
connecting adjacent vertices,
\be
  A_{x;i,j}
  = A_x (\xvec_{i,j} + \mbox{$\frac{1}{2}$} h_x \evec_x) , \quad
  A_{y;i,j}
  = A_y (\xvec_{i,j} + \mbox{$\frac{1}{2}$} h_y \evec_y) .
\ee
Here, $\evec_x$ and $\evec_y$ denote the unit vectors in the
$x$ and $y$ direction, respectively.
The definition of the discrete supercurrent $\Jvec_s$
is completely analogous.
The link variables, defined in Eq.~(\ref{Uvec}),
are obtained from the vector potential,
\be
U_{x;i,j} = {\rm e}^{ (\mathrm{i}/\kappa) A_{x;i,j} h_x } , \quad
U_{y;i,j} = {\rm e}^{ (\mathrm{i}/\kappa) A_{y;i,j} h_y } .
\ee
They are therefore also defined on the links.
Finally, the magnetic induction $\Bvec$,
which is a vector perpendicular to the plane
and given by the curl of the vector potential,
is defined at the center of a grid cell,
\be
  B_{z;i,j}
  = B_z (\xvec_{i,j} + \mbox{$\frac{1}{2}$} h_x \evec_x
  + \mbox{$\frac{1}{2}$} h_y \evec_y) .
\ee
The definition of the discrete variables
is also illustrated in Fig.~\ref{f-cell}.

Note that, because of the location of
the grid relative to the boundaries,
all scalar variables,
as well as the $x$ components of all vectors
($A_x$, $U_x$, $J_{s,x}$, and so forth),
are defined on a $N_x \times N_y$ grid,
whereas the $y$ components of all vectors
and the magnetic induction $B_z$
are defined on a $N_x \times (N_y -1)$ grid.

\paragraph{Boundary conditions.}
We assume periodicity in the $x$ direction,
so we need to consider the boundary conditions~(\ref{bc})
only at $y = y_1$ and $y = y_{N_y}$.

The boundary condition for the order parameter,
$\nvec \cdot \grad_{\bf A} \psi = 0$,
becomes
\be
U_{y;i,1} \psi_{i,2} - \psi_{i,1} = 0 , \quad
\psi_{i,N_y} - U_{y;i,N_y - 1}^* \psi_{i,N_y - 1} = 0 ,
\ee
for $i=1, \ldots\,, N_x$.
For the vector potential, we require that
$\partial_y A_x = H_z$ and $A_y$ is constant
($A_y = 0$) on the boundary.

\paragraph{Operators.}
The gradient of a scalar is a vector and
is therefore defined at the midpoint
of a link connecting two adjacent vertices.
Thus,
\be
  \left( \nabla \phi \right)_{x;i,j}
  = (\partial_x \phi) (\xvec_{i,j} + \mbox{$\frac{1}{2}$} h_x \evec_x)
  = h_x^{-1} (\phi_{i+1,j} - \phi_{i,j}) ,
\ee
with an analogous definition for the $y$ component.
The gauge-invariant ${\bf A}$-gradient
$\grad_{\bf A} = \grad + \mathrm{i}\Avec$
is defined in a similar way, with
\be
  \left( \grad_{\bf A} \psi \right)_{x;i,j}
  = h_x^{-1} (\psi_{i+1,j} U_{x;i,j} - \psi_{i,j}) .
\ee
Thus, the discrete version of the twisted
Laplacian~$\Delta_{\bf A}$ is
\begin{eqnarray}
  \left( \Delta_{\bf A} \psi \right)_{i,j}
  & = &
  h_x^{-2}
  (\psi_{i+1,j} U_{x;i,j} - 2\psi_{i,j} + \psi_{i-1,j} U_{x;i-1,j}^*)
  \nonumber \\
  &&+
  h_y^{-2}
  (\psi_{i,j+1} U_{y;i,j} - 2\psi_{i,j} + \psi_{i,j-1} U_{x;i,j-1}^*) .
\end{eqnarray}
The discrete version of the (normal) Laplacian
is defined in the usual way,
\be
  \left( \Delta \psi \right)_{i,j}
  = 
  h_x^{-2}
  (\psi_{i+1,j} - 2\psi_{i,j} + \psi_{i-1,j})
  +
  h_y^{-2}
  (\psi_{i,j+1} - 2\psi_{i,j} + \psi_{i,j-1}) .
\ee
The magnetic induction, which is the curl of
the vector potential, takes the form
\be
  B_{z;i,j}
  = h_x^{-1} (A_{y;i+1,j} - A_{y;i,j})
  - h_y^{-1} (A_{x;i,j+1} - A_{x;i,j}) .
\ee
We also need the divergence of the vector potential,
which is given by
\be
  \left( \div \Avec \right)_{i,j}
  = h_x^{-1} (A_{x;i,j} - A_{x;i-1,j})
  + h_y^{-1} (A_{y;i,j} - A_{y;i,j-1}) .
\ee

\paragraph{Algorithm.}
For numerical purposes, it is useful to treat
the TDGL equations (\ref{p}) and (\ref{A})
as two separate equations, which are coupled only
through certain fields and variables.
The electromagnetic potentials $\phi$ and $\Avec$
are treated as static variables in the order
parameter equation,
\begin{equation}
  \left( \partial_t - (\mathrm{i}/\kappa) \phi \right) \psi
  - \Delta_{\bf A} \psi - (1 - |\psi|^2) \psi
  = 0 .
\end{equation}
The local nonlinear part of this equation,
\begin{equation}
  \left( \partial_t - (\mathrm{i}/\kappa) \phi \right)\psi
  - (1 - |\psi|^2) \psi
  = 0 ,
\end{equation}
is integrated in the simplest possible manner,
\begin{equation}
  \psi_{i,j}(t + \Delta t)
  = {\rm e}^{-(i/\kappa) \phi_{i,j} \Delta t}
  \left\{ \psi_{i,j}(t) + \Delta t \left( 1 - |\psi_{i,j}|^2 \right) \psi_{i,j} \right\} .
\end{equation}
The nonlocal part,
\begin{equation}
\partial_t\psi - \Delta_{\bf A} \psi = 0 ,
\end{equation}
is integrated by using a backward Euler method,
where the linear equation system
is solved with the conjugate gradient method.

The equation for the vector potential,
\begin{equation}
  \sigma \partial_t \Avec
  - \Delta \Avec
  - (1/\kappa) \Jvec_s
  - \curl \Hvec
  = \mathbf{0} ,
\end{equation}
is linear and depends only indirectly on
the order parameter through the supercurrent.
If we treat the supercurrent as a static variable,
we can integrate the equation easily,
again using the backward Euler method.
In the actual implementation,
we use the fact that
the domain is periodic
to do a fast Fourier transform
in the $x$ direction, which
leaves us with a tridiagonal system
to solve in the $y$ direction.
This procedure is considerably faster
than using an iterative method,
such as the conjugate gradient method.

\subsection{Numerical results\label{ss-numerics}}
We use a rectangular sample,
periodic in the $x$ direction
and open in the $y$ direction,
with $N_x = N_y = 128$.
We take $h_x = h_y = \half \xi$,
so the sample measures 64 coherence lengths
in the periodic direction and
63.5 coherence lengths across.
(The coherence length $\xi$ is defined in
Eq.~(\ref{xi-lambda}).)

\begin{figure}[ht]
\begin{center}
\resizebox{4.0in}{!}{\mbox{
\includegraphics{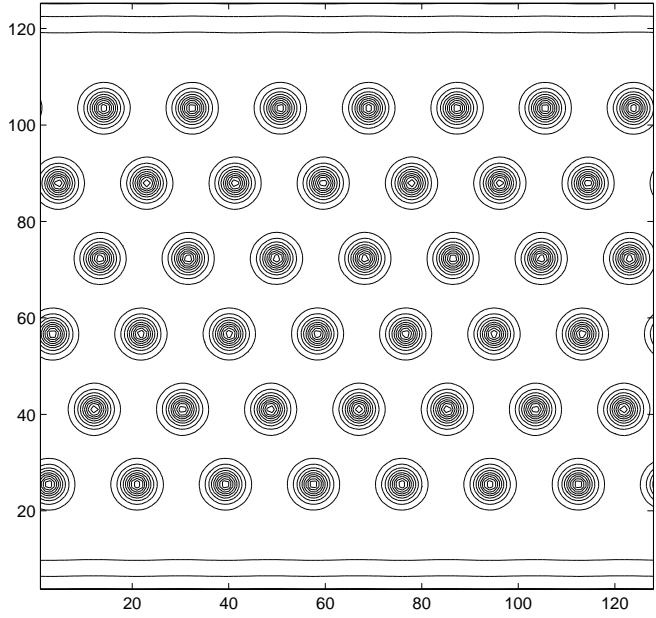}
}}
\caption{Contours of the density of Cooper pairs, $|\psi|^2$,
         for a system with $\kappa = 200$. \label{f-density}}
\end{center}
\end{figure}

First, we considered this system with $\kappa = 200$
and an applied magnetic field $H_z = 0.088 \kappa$.
With a relatively large value of $\kappa$,
the surface barrier for vortex entry is low,
and the system equilibrates relatively fast~\cite{ginzburg,chapman}.
The equilibration required $5\times10^4$ time steps
with $\Delta t = 0.4$.
The magnetic field produces an almost perfect vortex lattice.
Figure~\ref{f-density} gives a contour plot of
the density of Cooper pairs $|\psi|^2$ at equilibrium;
the zeros correspond to the centers of the vortices.

We then started from the configuration
of Fig.~\ref{f-density} to find
equilibrium configurations for
other values of $\kappa$,
varying $\kappa$ from $\kappa_{\min} = 40$
to $\kappa_{\max} = 800$.
In this range, the ground states are comparable 
and similar to the one shown in Fig.~\ref{f-density}.
Since the magnetization of a sample is proportional
to $1/\kappa^2$, the vortex density decreases with $\kappa$;
below $\kappa_{\min}$, the equilibrium state has fewer
vortices, and a comparison becomes meaningless.
Each equilibration required another $3\times10^4$ time steps.

\begin{figure}[ht]
\begin{center}
\resizebox{4.0in}{!}{\mbox{
\includegraphics{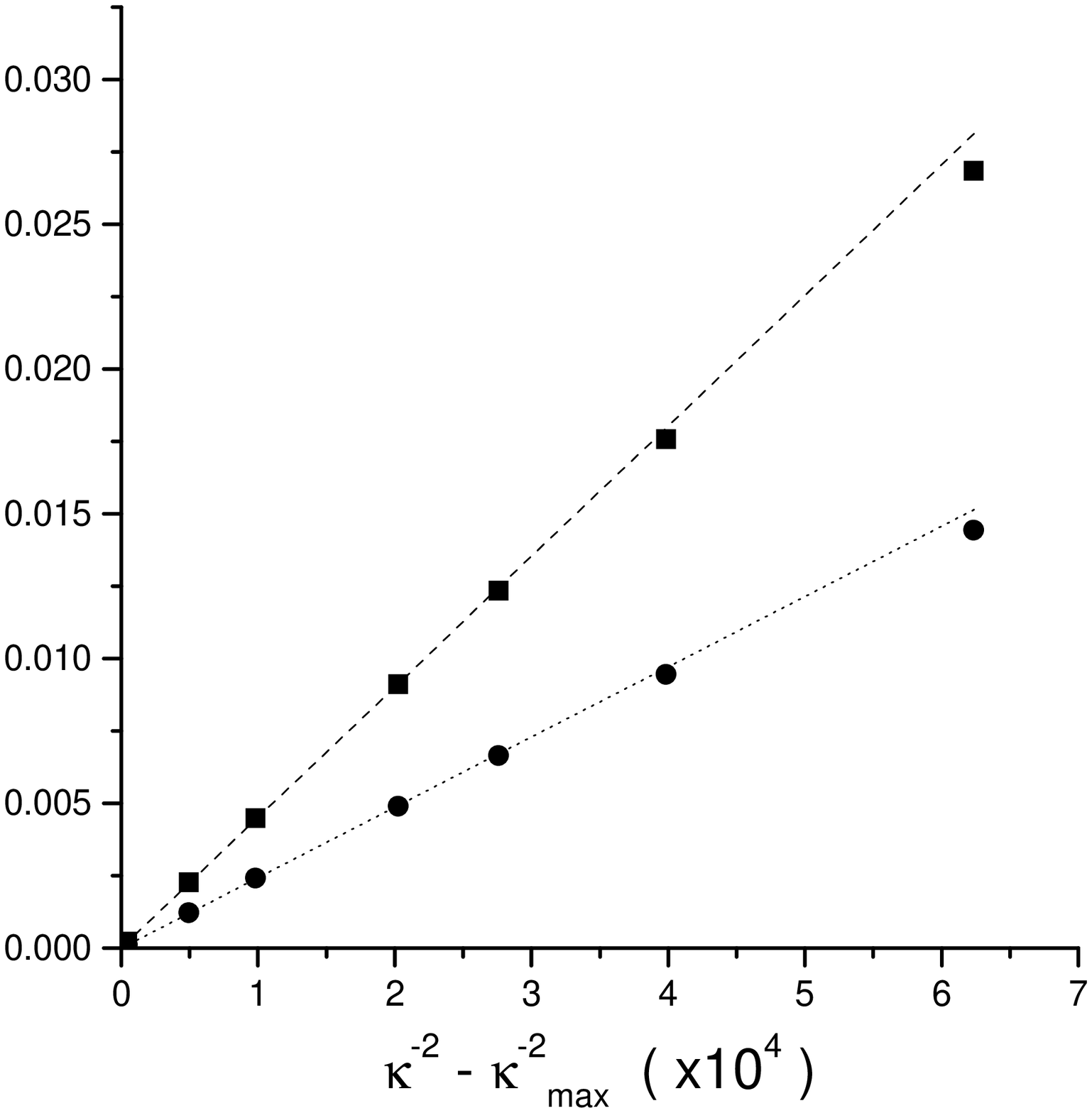}
}}
\caption{The quantities $\delta\psi$ (solid squares)
         and $\delta B_z$ (solid discs)
         for $\kappa = 40, 50, 60, 70, 100, 140, 400, 800$.
         The straight lines correspond to $1/\kappa^2$ behavior.
\label{f-asymp}}
\end{center}
\end{figure}

Figure~\ref{f-asymp}
gives the computed values
of the quantities
\be
  \delta \psi
  = \| \psi_\kappa  -  \psi_{\kappa_{\max}} \|_{L^2} , \quad
  \delta B_z
  = \frac{\| B_{z,\kappa} - B_{z, \kappa_{\max}} \|_{L^2}}
         {\| H_z \|_{L^2}} ,
\ee
for different values of $\kappa$.
The data show a behavior like $1/\kappa^2$
down to $\kappa \approx 40$.

\begin{figure}[ht]
\begin{center}
\resizebox{4.0in}{!}{\mbox{
\includegraphics{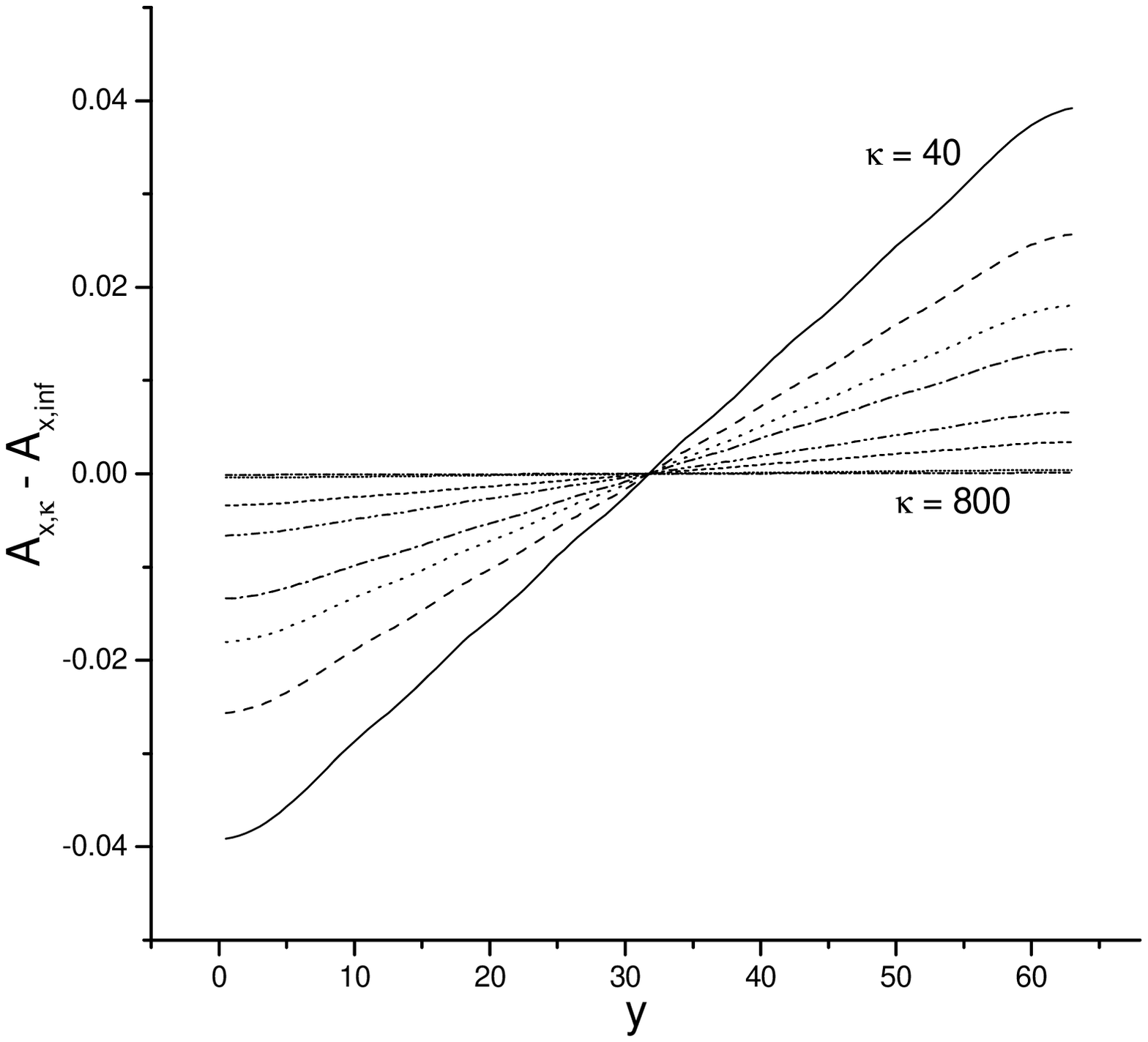}
}}
\caption{The average $<\!A_{x,\kappa} - A_{x,\kappa_{\max}}\!>$ vs.\ $y$
for $\kappa = 40, 50, 60, 70, 100, 140, 400, 800$.
\label{f-Aavg}}
\end{center}
\end{figure}

Figure~\ref{f-Aavg}
shows the average over $x$ of
$A_{x,\kappa} - A_{x,\kappa_{\max}}$
as a function of $y$ in the bulk of the sample,
for different values of $\kappa$.

The numerical results show that the solution
of the TDGL equations converges as $\kappa$
increases; in fact, they show quadratic convergence
in the small parameter $1/\kappa$.
Given the fact that the Ginzburg--Landau parameter
of high-$T_c$ superconducting materials is
of the order of 50--100,
we conclude that the limiting equation is
a practical alternative in many applications.
The question thus becomes:
What is the limiting equation,
and can we confirm the numerical conclusions
by rigorous arguments?
We address this question in the next section.

\section{Asymptotic analysis\label{s-aa}}
\setcounter{equation}{0}
We now return to the TDGL equations~(\ref{tdgl-p})--(\ref{tdgl-bc})
and consider their limit as $\kappa \to \infty$.
These are our standing hypotheses:
\begin{description}
\item[(H1)] $\Omega$ is bounded in ${\bf R}^n$ ($n =2,3$),
with a sufficiently smooth boundary $\partial \Omega$,
for example, $\partial \Omega$ of class $C^{1,1}$.
\item[(H2)] The parameters $\kappa$ and $\sigma$
are real and positive.
\item[(H3)] $\Hvec$ is independent of time;
as a function of position, it satisfies the regularity
condition $\Hvec \in [W^{\alpha,2}(\Omega)]^n$
for some $\alpha \in (\half, 1)$.
\item[(H4)] $\kappa \gg 1$;
$\sigma = O(1)$ and
$\Hvec = O(\kappa)$
as $\kappa \to \infty$.
\end{description}
The assumptions~\textbf{(H1)--(H3)} suffice
to prove that the TDGL equations define a
dynamical system in the Hilbert space
\be
  \mathcal{W}^{1+\alpha,2}
  =
  [W^{1+\alpha,2} (\Omega)]^2 \times [W^{1+\alpha,2} (\Omega)]^n ;
\ee
see~\cite{Fleck-Kaper-Takac-98}.
The space $W^{1+\alpha,2} (\Omega)$
is continuously imbedded in
$W^{1,2} (\Omega) \cap L^\infty (\Omega)$,
so $\psi$ and $\Avec$ are bounded and differentiable
with square-integrable (generalized) derivatives.
\textbf{(H4)} is the operative hypothesis
for the asymptotic analysis.

\subsection{Mathematical analysis \label{ss-analysis}}
\paragraph{Scaling.}
We begin by scaling the TDGL equations,
taking into account the fact that
we are interested in the limit
as $q_s \to 0$ (\textit{weak coupling}),
when the applied field is near
the upper critical field.
The scaling is done by means of
the dimensionless GL parameter $\kappa$,
which grows like $|q_s|^{-1}$.

Since $\Hvec = O(\kappa)$ as $\kappa \to \infty$,
we begin by scaling $\Hvec$ by a factor $\kappa$,
$\Hvec = \kappa \Hvec'$.
By scaling the vector potential by
the same factor $\kappa$,
we achieve that the electromagnetic variables
are all of the same order.

The scalar potential is proportional to
the charge density of the Cooper pairs,
which is $O(|q_s|)$ as $q_s \to 0$.
Hence, $\kappa \phi$ remains of order one.
This suggests scaling $\phi$ by a factor
$\kappa^{-1}$.

Table~\ref{scale} summarizes the relation
between the current (nondimensional) variables
and their scaled (primed) counterparts.
We adopt the latter as the new variables
and work until further notice
on the scaled problem.
We omit all primes.
\begin{table}[ht]
\begin{center}
\begin{footnotesize}
\caption{Scaling.
  \label{scale}}\vspace{1ex}
\begin{tabular}{||l|l||}\hline
Applied Field
  & $\Hvec = \kappa \Hvec'$ \\\hline
  & $\psi = \psi'$ \\
Dependent variables
  & $\Avec = \kappa \Avec'$ \\
  & $\phi = \kappa^{-1} \phi'$ \\\hline
  & $\Bvec = \kappa \Bvec'$ \\
Electromagnetic variables
  & $\Jvec = \kappa \Jvec'$ \\
  & $\Evec = \kappa \Evec'$ \\\hline
\end{tabular}
\end{footnotesize}
\end{center}
\end{table}

After scaling, the relevant parameter is $\kappa^2$,
rather than $\kappa$, so we introduce $\eps$,
\be
  \eps = \kappa^{-2} .
  \Label{eps}
\ee
The scaled TDGL equations are
\be
  \left( \partial_t - \mathrm{i} \eps \phi \right) \psi
  - (\grad + \mathrm{i} \Avec)^2 \psi
  - (1 - | \psi |^2) \psi
  = 0 ,
  \Label{p}
\ee
\be
  \sigma \partial_t \Avec
  - \Delta \Avec
  - \eps \Jvec_s
  - \curl \Hvec
  = {\bf 0} ,
  \Label{A}
\ee
where
\be
  \Jvec_s
  =
  - \frac{1}{2\mathrm{i}}
  (\psi^* \grad \psi - \psi \grad \psi^*)
  -
  |\psi|^2 \Avec
  =
  - \Im
  \left[ \psi^* (\grad + \mathrm{i} \Avec) \psi \right] ,
  \Label{Js}
\ee
with the corresponding gauge condition,
\be
  \eps \sigma \phi +  \div \Avec = 0 .
 \Label{gauge}
\ee
The boundary conditions associated with
Eqs.~(\ref{p}) and (\ref{A}) are
\be
  \nvec \cdot (\grad + \mathrm{i} \Avec) \psi = 0 , \quad
  \nvec \times (\curl \Avec) = \nvec \times \Hvec , \quad
  \nvec \cdot \Avec = 0 .
  \Label{bc}
\ee
The electromagnetic variables are given by the expressions
\be
  \Bvec = \curl \Avec , \quad
  \Jvec = \curl\curl \Avec , \quad
  \Evec = - \partial_t \Avec - \eps \grad \phi .
  \Label{BJE}
\ee

\paragraph{Reduction to homogeneous form.}
Next, we homogenize the problem.
Let $\Avec_0$ be the (unique) minimizer
of the convex quadratic form $Q_1 \equiv Q_1 [\Avec]$,
\be
  Q_1 [\Avec]
  = 
  \int_\Omega
  \left[ (\div \Avec)^2 + |\curl \Avec - \Hvec|^2 \right]
  \ {\rm d}{\xvec} ,
  \Label{QAH}
\ee
on dom$(Q_1) = \{ \Avec \in [W^{1,2}(\Omega)]^n :
\nvec \cdot \Avec = 0
\mbox{ on } \partial \Omega \}$.
This minimizer satifies the boundary-value problem
\be
  \curl \curl \Avec
  - \curl \Hvec
  = {\bf 0} , \quad
  \div \Avec = 0 \quad
  \mbox{in } \Omega ,
  \Label{eq-A0}
\ee
\be
 \nvec \times (\curl \Avec) = \nvec \times \Hvec , \quad
 \nvec \cdot \Avec = 0 \quad
  \mbox{on } \partial \Omega ,
  \Label{bc-A0}
\ee
in the dual of dom$(Q_1)$ with respect to
the inner product in $[L^2(\Omega)]^n$.
The mapping $\Hvec \mapsto \Avec_0$
is linear, time independent, and continuous
from $[W^{\alpha,2}(\Omega)]^n$
to $[W^{1+\alpha,2}(\Omega)]^n$~\cite{Georgescu}.
The contribution of the vector $\Avec_0$
to the magnetic field is
\be
  \Bvec_0 = \curl \Avec_0 .
  \Label{B0}
\ee
We substitute variables,
\be
  \Avec = \Avec_0 + \eps \Avec' ,
  \Label{A'}
\ee
and rewrite the (scaled) TDGL equations~(\ref{p})--(\ref{bc})
in terms of $\psi$ and $\Avec'$ (omitting the primes),
\be
  \partial_t \psi + \mathrm{i} \sigma^{-1} (\div (\eps \Avec)) \psi
  - (\grad + \mathrm{i} (\Avec_0 + \eps \Avec))^2 \psi
  - (1 - | \psi |^2) \psi
  = 0
  \quad \mbox{in } \Omega ,
  \Label{p-red}
\ee
\be
  \sigma \partial_t \Avec
  - \Delta \Avec
  - \Jvec_s
  = {\bf 0}
  \quad \mbox{in } \Omega ,
  \Label{A-red}
\ee
where
\be
  \Jvec_s
  =
  - \frac{1}{2\mathrm{i}} (\psi^* \grad \psi - \psi \grad \psi^*)
  - |\psi|^2 (\Avec_0 + \eps \Avec) ,
  \Label{Js-red}
\ee
and
\be
  \nvec \cdot \grad \psi = 0 , \quad
  \nvec \times (\curl \Avec) = {\bf 0} , \quad
  \nvec \cdot \Avec = 0
  \quad \mbox{on } \partial \Omega .
  \Label{bc-red}
\ee

\paragraph{Functional formulation.}
We reformulate the system of
Eqs.~(\ref{p-red})--(\ref{bc-red})
as an ordinary differential equation
for a vector-valued function
$u = (\psi, \Avec)$
from the time domain $(0, \infty)$
to a space of functions on $\Omega$,
\be
  u = (\psi, \Avec) :
  [0,\infty) \to \mathcal{L}^2
  = [L^2 (\Omega)]^2 \times [L^2 (\Omega)]^n .
\ee
The equation is
\be
  \frac{{\rm d}u}{{\rm d}t} + Au = f_0(u) + \eps f_1 (u) ,
  \Label{eq-ueps}
\ee
where $A$ is the linear operator in $\mathcal{L}^2$
associated with the quadratic form $Q \equiv Q[u]$,
\be
  Q[u] = \int_\Omega
  \left[
  |\grad \psi|^2 + \sigma^{-1} \left( (\div \Avec)^2 + |\curl \Avec|^2 \right)
  \right] \ {\rm d}\xvec ,
  \Label{Qu}
\ee
on dom$(Q) = \{ u = (\psi, \Avec) \in \mathcal{L}^2 :
\nvec \cdot \Avec = 0 \mbox{ on } \partial \Omega \}$.
The functions $f_0$ and $f_1$ are nonlinear,
\be
  f_i(u) = (\varphi_i (\psi, \Avec), \sigma^{-1} \Fvec_i (\psi, \Avec)) , \quad
  i=0,1 ,
\ee
where
\begin{eqnarray}
  \varphi_0 (\psi, \Avec)
  &=& 2\mathrm{i} \Avec_0 \cdot (\grad \psi)
  - | \Avec_0 |^2 \psi
  + (1 - |\psi|^2) \psi , \\
  \varphi_1 (\psi, \Avec)
  &=& \mathrm{i}(1 - \sigma^{-1}) (\div \Avec) \psi
  + 2\mathrm{i} \Avec \cdot (\grad \psi)
  - (\Avec_0 \cdot \Avec) \psi - |\Avec|^2 \psi , \\
  \Fvec_0 (\psi, \Avec)
  &=& 0 , \\
  \Fvec_1 (\psi, \Avec)
  &=&
  - \frac{1}{2\mathrm{i}}
  (\psi^* \grad \psi - \psi \grad \psi^*)
  -
  |\psi|^2 (\Avec_0 + \eps \Avec) .
\end{eqnarray}
Given any $f = (\varphi, \sigma^{-1} \Fvec) \in \mathcal{L}^2$,
the equation $Au = f$ is equivalent with
the system of uncoupled boundary-value problems
\begin{eqnarray}
  - \Delta \psi = \varphi \;
  \mbox{ in } \Omega , &&
  \nvec \cdot \grad \psi = 0 \;
  \mbox{ on } \partial \Omega , \\
  - \Delta \Avec = \Fvec \;
  \mbox{ in } \Omega , &&
  \nvec \times \Avec = {\bf 0} , \;
  \nvec \cdot \Avec = 0 \;
  \mbox{ on } \partial \Omega ,
\end{eqnarray}
in the dual of dom$(Q)$ with respect to
the inner product in $\mathcal{L}^2$.
The operator $A$ is selfadjoint and
positive definite in $\mathcal{L}^2$;
hence, its fractional powers $A^{\theta/2}$
are well defined, they are unbounded
if $\theta \geq 0$, and
dom$(A^{\theta/2})$ is a closed linear subspace
of $\mathcal{W}^{\theta,2}
= [W^{\theta,2} (\Omega)]^2 \times [W^{\theta,2} (\Omega)]^n$;
see~\cite[Section~1.4]{Henry}.

The solution of Eq.~(\ref{eq-ueps}) depends on $\eps$;
we denote it by $u_\eps$.
We compare $u_\eps$ with the solution $u_0$ of the reduced equation
\be
  \frac{{\rm d}u}{{\rm d}t} + Au = f_0(u) .
  \Label{eq-u0}
\ee

\begin{THEOREM}  \Label{conv}
There exists a positive constant $C$ such that
\be
  \| u_\eps (t) - u_0 (t) \|_{\mathcal{W}^{1+\alpha,2}}
  \leq
  C \left(
  \| u_\eps (0) - u_0 (0) \|_{\mathcal{W}^{1+\alpha,2}}
  + \eps \right) ,
  \quad t \in [0, T] .
\ee
\end{THEOREM}

\begin{PROOF}
Let $B_R$ be the ball of radius $R$ centered at the origin
in $\mathcal{W}^{1+\alpha,2}$.
Let $u_\eps \in B_R$ and $u_0 \in B_R$ satisfy
Eqs.~(\ref{eq-ueps}) and (\ref{eq-u0}), respectively,
with initial data $u_\eps (0)$ and $u_0 (0)$.
The difference $v = u_\eps - u_0$
satisfies the differential equation
\be
  \frac{{\rm d}v}{{\rm d}t} + Av = f_0(u_\eps) - f_0(u_0) + \eps f_1(u_\eps)
  \Label{eq-v}
\ee
or, equivalently, the integral equation
\be
  v(t) = {\rm e}^{-tA} v(0)
  + \int_0^t {\rm e}^{-(t-s)A} [f_0(u_\eps) - f_0(u_0) + \eps f_1(u_\eps)](s)\ {\rm d}s .
  \Label{v-int}
\ee
From the integral equation we obtain the estimate
\begin{eqnarray}
  \| v(t) \|_{\mathcal{W}^{1+\alpha,2}} &\leq&
  \| {\rm e}^{-tA} \|_{\mathcal{W}^{1+\alpha,2}}
  \| v(0) \|_{\mathcal{W}^{1+\alpha,2}}
  + \int_0^t \|A^{(1+\alpha)/2} {\rm e}^{-(t-s)A} \|_{\mathcal{W}^{1+\alpha,2}}
  \nonumber\\
  &&\times \left[ \|f_0(u_\eps) - f_0(u_0)\|_{L^2}
  + \eps \|f_1(u_\eps)\|_{L^2}\right](s)\ {\rm d}s .
  \Label{v-est}
\end{eqnarray}
The operator norms satisfy the inequalities
\be
  \| {\rm e}^{-tA} \|_{\mathcal{W}^{1+\alpha,2}}
  \leq 1 , \quad
  \| A^{(1+\alpha)/2} {\rm e}^{-(t-s)A} \|_{\mathcal{W}^{1+\alpha,2}}
  \leq
  C (t-s)^{-(1+\alpha)/2} ;
  \Label{ops-est}
\ee
see~\cite[Theorem 1.4.3]{Henry}.
Furthermore, adding and subtracting terms, we have
\begin{eqnarray}
  f_0(u_\eps) - f_0(u_0)
  &=&
  \left(
  2i \Avec_0 \cdot (\grad(\psi_\eps - \psi_0))
  - |\Avec_0|^2 (\psi_\eps - \psi_0) \right. \nonumber \\
  &&+\left. (1 - |\psi_\eps|^2 - |\psi_0|^2) (\psi_\eps - \psi_0)
   - \psi_\eps \psi_0 (\psi_\eps^* - \psi_0^*) , \
  0
  \right) ,
\end{eqnarray}
where
\[
  \| 2i \Avec_0 \cdot (\grad(\psi_\eps - \psi_0)) \|_{L^2}
  \leq
  2 \| \Avec_0 \|_{L^\infty} 
  \| \psi_\eps - \psi_0 \|_{W^{1,2}}
\]
\[
  \leq
  C \| \psi_\eps - \psi_0 \|_{W^{1+\alpha,2}}
  \leq
  C \| u_\eps - u_0 \|_{\mathcal{W}^{1+\alpha,2}} ,
\]
\[
  \| |\Avec_0|^2 (\psi_\eps - \psi_0) \|_{L^2}
  \leq
  C \| \Avec_0 \|_{L^\infty}^2
  \| \psi_\eps - \psi_0 \|_{L^\infty}
\]
\[
  \leq
  C \| \psi_\eps - \psi_0 \|_{W^{1+\alpha,2}}
  \leq
  C \| u_\eps - u_0 \|_{\mathcal{W}^{1+\alpha,2}} ,
\]
and the other terms are estimated similarly.
Here, $C$ is some (generic) positive constant,
which may depend on $\Hvec$ and $\Omega$ but not
on $u_\eps$ or $u_0$.
(In these inequalities we have used
the continuity of the imbedding of
$W^{1+\alpha,2}$ into $W^{1,2} \cap L^\infty$.)
The result is an inequality of the type
\be
  \| f_0(u_\eps) - f_0(u_0) \|_{L^2}
  \leq
  C \| u_\eps - u_0 \|_{\mathcal{W}^{1+\alpha,2}} ,
  \Label{f0-Lip}
\ee
showing that $f_0$ is Lipschitz from
${\mathcal{W}^{1+\alpha,2}}$ to $\mathcal{L}^2$.

Using similar estimates, we show that
$f_1$ is bounded from ${\mathcal{W}^{1+\alpha,2}}$
to $\mathcal{L}^2$, so there exists a positive
constant $C$ such that
\be
  \| f_1 (u_\eps) \|_{L^2} \leq C .
  \Label{f1-bd}
\ee
Combining the estimates~(\ref{ops-est}), (\ref{f0-Lip}),
and (\ref{f1-bd}) with the inequality~(\ref{v-est}),
we conclude that there exist positive constants $C_1$ and $C_2$
such that
\be
  \| v(t) \|_{\mathcal{W}^{1+\alpha,2}}
  \leq
  \| v(0) \|_{\mathcal{W}^{1+\alpha,2}}
  + \eps C_1 t^{(1-\alpha)/2}
  + C_2
  \int_0^t (t-s)^{-(1+\alpha)/2}
  \| v(s) \|_{\mathcal{W}^{1+\alpha,2}}
  \ {\rm d}s .
\ee
Applying Gronwall's inequality, we obtain the estimate
\be
  \| v(t) \|_{\mathcal{W}^{1+\alpha,2}}
  \leq
  C
  \left(\| v(0) \|_{\mathcal{W}^{1+\alpha,2}} + \eps \right) , \quad
  t \in [0,T] ,
\ee
for some positive constant $C$.
\end{PROOF}

It follows from Theorem~\ref{conv} that,
if the initial data are such that
$\| u_\eps (0) - u_0 (0) \|_{\mathcal{W}^{1+\alpha,2}} = o(1)$
as $\eps \downarrow 0$, then
\be
  \lim_{\eps \to 0} u_\eps  = u_0
  \Label{u-conv}
\ee
in $C([0,T]; \mathcal{W}^{1+\alpha,2})$ for any $T>0$.
In particular, if
$\| u_\eps (0) - u_0 (0) \|_{\mathcal{W}^{1+\alpha,2}} = O(\eps)$,
then the convergence in Eq.~(\ref{u-conv}) is $O(\eps)$.

\subsection{Interpretation and final remarks\label{ss-interp}}
It remains to translate the results back
in terms of the original variables.
We denote the solution of the TDGL equations,
Eqs.~(\ref{tdgl-p})--(\ref{tdgl-bc}),
by $\psi_\kappa$, $\Avec_\kappa$, $\phi_\kappa$.
The variables $\Avec_\kappa$ and $\phi_\kappa$
are related by the gauge condition
$\sigma \phi_\kappa + \div \Avec_\kappa = 0$
at all times.
Let $\Bvec_\kappa = \curl \Avec_\kappa$.

Let $\Avec_\infty$ be the solution of
the boundary-value problem
\be
  \curl \curl \Avec
  - \curl \Hvec
  = {\bf 0} , \quad
  \div \Avec = 0 \quad
  \mbox{in } \Omega ,
  \Label{eq-Ainf}
\ee
\be
 \nvec \times (\curl \Avec) = \nvec \times \Hvec , \quad
 \nvec \cdot \Avec = 0 \quad
  \mbox{on } \partial \Omega ,
  \Label{bc-Ainf}
\ee
and put $\Bvec_\infty = \curl \Avec_\infty$.
The vector $\Avec_\infty$ and, hence, $\Bvec_\infty$
do not vary with time.
Let $\psi_\infty$ satisfy the equations
\be
  \partial_t \psi
  - \Delta_{{\bf A}_\infty} \psi
  - (1 - | \psi |^2) \psi
  = 0 \;
  \mbox{ in } \Omega , \quad
  \nvec \cdot \grad_{{\bf A}_0} \psi
  = 0 \;
  \mbox{ on } \partial\Omega .
  \Label{pinf}
\ee
Then it follows from Theorem~\ref{conv} that
there exists a positive constant $C$ such that
\begin{eqnarray}
\lefteqn{
  \| \psi_\kappa (t) - \psi_\infty (t) \|_{W^{1+\alpha,2}}
  +
  \frac{\| \Bvec_\kappa (t) - \Bvec_\infty \|_{W^{\alpha,2}}}
       {\| \Hvec \|_{W^{\alpha,2}}}
  } \\
  &&\leq
  C \left(
  \| \psi_\kappa (0) - \psi_\infty (0) \|_{W^{1+\alpha,2}}
  +
  \frac{\| \Bvec_\kappa (0) - \Bvec_\infty \|_{W^{\alpha,2}}}
       {\| \Hvec \|_{W^{\alpha,2}}}
  +
  \frac{1}{\kappa^2}
  \right) ,
\end{eqnarray}
for all $t \in [0,T]$, $T > 0$.

The approximation $(\psi_\infty, \Bvec_\infty)$
is the ``frozen-field approximation.''
Hence, the analysis shows that the solution
of the TDGL equations converges to
the frozen-field approximation,
uniformly on compact time intervals $[0, T]$
in the topology of $[W^{1+\alpha,2} (\Omega)]^2 \times
[W^{\alpha,2} (\Omega)]^n$,
as soon as the initial data satisfy the asymptotic estimates
$\| \psi_\kappa (0) - \psi_\infty (0) \|_{W^{1+\alpha,2}} = o(1)$
and
$\| \Bvec_\kappa (0) - \Bvec_\infty \|_{W^{\alpha,2}} = o(\kappa)$
as $\kappa \to \infty$.
Under slightly sharper conditions
we obtain the order of convergence.

\begin{COROLLARY}
\label{corr}
If
\[
  \| \psi_\kappa (0) - \psi_\infty (0) \|_{W^{1+\alpha,2}}
  = O \left( \frac{1}{\kappa^2} \right)
  \quad\mbox{and}\quad
  \frac{\| \Bvec_\kappa (0) - \Bvec_\infty \|_{W^{\alpha,2}}}
       {\| \Hvec \|_{W^{\alpha,2}}}
  = O \left( \frac{1}{\kappa^2} \right)
\]
as $\kappa \to \infty$,
then
\be
  \| \psi_\kappa (t) - \psi_\infty (t) \|_{W^{1+\alpha,2}}
  +
  \frac{\| \Bvec_\kappa (t) - \Bvec_\infty \|_{W^{\alpha,2}}}
       {\| \Hvec \|_{W^{\alpha,2}}}
  =
  O \left( \frac{1}{\kappa^2} \right) ,
  \label{AH-conv}
\ee
uniformly on compact intervals.
\end{COROLLARY}

This result explains the numerical results
presented in Section~\ref{s-num}.

\paragraph{Remark~1.}
The asymptotic approximation procedure can be continued
to higher order, as can be seen from a formal expansion.
The equations for the order parameter and the vector potential
decouple, and at each order one finds first the vector potential,
then the order parameter.
The vector potential satisfies a linear heat equation;
for example, the first correction beyond $\Avec_\infty$
is $\kappa^{-1} \Avec$, where $\Avec$ satisfies the equation
\be
  - \sigma \partial_t \Avec
  + \Delta \Avec
  = \Im \left[ \psi_\infty^* \grad_{{\bf A}_\infty} \psi_\infty \right] .
\ee
\paragraph{Remark~2.}
The analysis given here differs at several points
from the analysis of Ref.~\cite{Du-Gray-96}.
First, our scaling is slightly different and,
we believe, more in tune with the physics;
second, our regularity assumptions on the
applied field are weaker;
third, our proofs are more direct;
and fourth, our results hold in a stronger topology.

\subsection*{Acknowledgments}
We thank Professor Todd Dupont (University of Chicago)
for stimulating discussions throughout the course of this investigation.
We also acknowledge the work of Damien Declat (student
from ISTIL, Lyon, France), who assisted in the development
of an earlier version of the parallel computer program.

This work was supported by the Mathematical, Information,
and Computational Sciences Division subprogram of the
Office of Advanced Scientific Computing Research, U.S. Department
of Energy, under Contract W-31-109-Eng-38.
The second author was partially supported by
the University of Chicago/Argonne National Laboratory
Collaborative Grant No.~96-011.

\end{document}